# Chaos and chaos control in MEMS resonators under power law noise

Yunfang Wang[1] and Youming Lei[1, 2*]

[1] Department of Applied Mathematics, Northwestern Polytechnical University, Xi'an 710072, China

[2] MIIT Key Laboratory of Dynamics and Control of Complex Systems, Northwestern Polytechnical University, Xi'an 710072, China

## Abstract

The chaotic dynamics and its control under power law noise in Micro-electromechanical Systems (MEMS) resonators with electrostatic excitation are probed. On the basis of the stochastic Melnikov method in the mean-square sense and the mean largest Lyapunov exponent, the threshold value of power law noise intensity for the onset of chaos is obtained analytically and numerically. We show that the threshold of noise intensity decreases with the increasing of the frequency exponent of power law noise in the parameter space. Numerical simulations, such as phase diagram and time history, are employed to verify the results acquired by the stochastic Melnikov method. Inspired by the analytical results, a time-delay feedback control algorithm is proposed for controlling the chaotic motion in the resonators. The effectiveness of this controller is certificated by the above numerical method.

Keywords: MEMS resonator; Power law noise; Stochastic Melnikov method; Chaos control

## 1. Introduction

The Micro-electromechanical resonators have evolved into a serviceable instrument for Micro-electromechanical Systems (MEMS) to exchange messages and communicate with outside world [1]. The electrostatic actuation is a high-efficiency non-contact drive mode that can avoid lowering the mechanical properties of sensors, so it becomes the preferred drive mode for a mass of MEMS resonators. Based on this, the mathematical model of electrostatic driven micro-beam is introduced in subsection 2.1. For many MEMS, frequency control and frequency selection are the main functions and applications of resonators. Once the resonators are deemed as the frequency source by producing resonant frequency and conducting frequency control, the existence of noise will exert a significant impact on the properties and nature of MEMS, such as limiting the receiver's channel spacing and selectivity, causing bit errors in digital communications system, limiting the accuracy of synchronization and syntonization [2], and limiting the output at low values of input(control) signals [3]. On account of the miniaturization and integration of resonators, such problems can be especially acute. In view of the unavoidable inherent nature of noise on electronic devices, it is a judicious choice to study impacts of noise on the responses of MEMS resonators and how to ameliorate and eliminate them. Based on the same consideration, Jokić et al. [4] explored the dependence of Adsorption-Desorption Phase Noise on the resonant frequency, operating pressure, and temperature in MEMS/NEMS resonators. Zhang et al. [5] analyzed the necessary conditions for of chaos induced by bounded noise in MEMS resonators by the Melnikov method and verified the experimental results by bifurcation diagrams, Poincare maps, phase portraits and time histories. Handel [6] provided the resulting formula for physical $1/f$ noise in MEMS resonators and sensors, etc. But the type of noise mentioned above isn't ubiquitous in electronic devices compared with power law noise. Power law noise, generated by the random fluctuation of carrier density in the active devices, is the main noise type near the center frequency and in the range of low frequency with large amplitude, so it must be taken into account when designing device models. Furthermore, power law noise maybe both the most interesting but vexing of all observed noise [7]. Although the noise was detected decades ago, there is no general theory to explain it, including its model and origin. Hence it continues to attract significant interest through abundant theoretical and experimental work to understand the mechanism behind it. Chaotic behaviors induced by noise are a class of nonlinear motion, which is notorious in the dynamic response of numerous systems [8]; however it has been verified to be practical in plentiful of areas, sensing, fluid mixture, and secure communication [9,10]. For the purpose of enhancing the work efficiency of MEMS resonators, inducing or restraining the generation of chaos, the chaotic motion, chaos areas and chaos control measures become a pressing issue. Diverse chaotic behaviors of MEMS resonators caused by power law noise and chaos control method will be set to the focus of this work.

[*] Corresponding author. E-mail Address: leiyouming@nwpu.edu.cn

With regard to the nonlinear dynamic behaviors and chaos control methods of MEMS resonators, a great number of investigators and pursuers have made abundant pioneering heuristic efforts and attempts. Liu, Davidson and Lin [11] investigated both the static and dynamic instabilities of MEMS cantilever subject to weak and strong disturbances. They observed period doubling, chaos, and strange attractors for both open and closed loop control systems with strong disturbances. Farid et al. [12] implemented the Melnikov method and the maximum velocity criterion to derive a necessary condition for the initiation of chaos in the electrostatically actuated arch micro-nano resonators. Bifurcation diagram, Lyapunov exponent and Poincare section confirm the validity of their analytical expressions. Albert and Wang [13] investigated the chaotic motion in a certain frequency band of simplified MEMS devices and simultaneously determined the corresponding equilibrium, natural frequency and response to provide predictions for design, manufacturing, testing and industrial applications. Ehsan et al. [14] employed diversified methods, such as the multiple scales method, the Melnikov method and a novel method proposed by themselves, to predict chaos in MEMS-NEMS resonators and they believed that the results obtained with the Melnikov method can be applied to design and acquire the optimum operational conditions. Ding et al. [15] studied the dynamics of a delayed MEMS nonlinear coupled system. By using the multi-time scale method, they derived the normal form near the Hopf and Hopf-pitchfork bifurcations critical points and showed the region near the bifurcation critical point where the MEMS nonlinear coupled system has a stable fixed point or a stable periodic solution. On the basis of chaos control tactics, the responses of electrostatically actuated MEMS resonators have been investigated far and wide. Haghighi and Markazi [16] utilized the Melnikov method to inquire the prediction of chaos in MEMS resonators under electrostatic force and applied the robust adaptive fuzzy control method to stabilize the MEMS beam. In the same way, Siewe and Hegazy [17] studied the analytic criterion for homoclinic chaos and introduced a time-varying stiffness to control the chaotic motion. Alhababa [18] investigated the existence of chaos in a non-autonomous fractional-order MEMS resonator using the maximal Lyapunov exponent and plotting the strange attractors. Afterward, the novel fractional finite-time controller was introduced to suppress the chaos arising from model uncertainties and external disturbance in a given finite time. Kwangho et al. [10] put forward that the nonlinear coupling between applied electrostatic force and the mechanical motion of the resonators can lead to chaotic oscillations and proposed a control strategy to convert chaos into periodic motion with enhanced output energy. Luo et al. [19] focused on the adaptive control of the arch MEMS resonator. In the controller design, what was novel and different from other researchers was the introduction of the Chebyshev neural network system to learn unknown dynamical behavior. By combining the observer, first-order filter, neural network and Nussbaum function, an adaptive control method was proposed to force the system state to appropriate the reference signal and suppress the chaotic oscillation of the arch MEMS resonator near the mean resonance frequency. In another paper, Luo et al. [20] proposed an observer-based adaptive stabilization scheme. The frequency distribution model of the fractional integrator and the fractional Lyapunov stability criterion were used to realize the stability of the fractional-order chaotic MEMS resonator under the condition of uncertain function, parameter perturbation and unmeasurable states with electrostatic excitation. Compared with integer-order MEMS resonators, fractional-order systems could simulate their genetic properties better and exhibit complex dynamic behaviors. Yau et al. [21] employed the phase portrait, maximum Lyapunov exponent and bifurcation diagram to find the chaotic areas. In order to suppress chaos, a robust fuzzy sliding mode controller is designed to turn the chaotic motion into a periodic motion. Tusset et al. [22] introduced three control strategies: State Dependent Riccati Equation (SDRE) Control, Optimal Linear Feedback Control and Fuzzy Sliding Mode Control for controlling the trajectories in a fractional order dynamic system.

Based on the extension of the stochastic Melnikov process, in this study, the stochastic Melnikov method will be the pivotal manner to study chaotic behaviors and chaos control in MEMS resonators under power law noise. The rest of the paper is organized as follows. In Section 2, we will draw forth the mathematical model of clamped-clamped microbeam and describe the power law noise in detail. Section 3 presents a threshold curve of power law noise intensity versus the frequency exponent with the stochastic Melnikov method. Numerical results will be utilized to verify the analytical findings and to investigate the effects of noise on the chaotic behaviors in Section 4. In Section 5, a time-delay feedback control is used to control chaotic motion. Summaries and conclusions are drawn in Section 6.

## 2. Mathematical model and power law noise

The mathematical model of investigated MEMS resonator with electrostatic actuation mode and stochastic characteristics of power law noise will be introduced briefly in this section.

### 2.1 Mathematical model

With regard to a class of clamped-clamped beam resonators, just as we all know, is actuated by an external driving force with direct current bias voltage between the electrode and the resonator. In this work, we elect the dynamic model of MEMS resonator derived by Haghighi and Markazi [16] for further research. The governing equation of motion can be expressed as

$$m\ddot{z} + b\dot{z} + k_1 z + k_3 z^3 = F(z, \tau) + G(\tau), \qquad (1)$$

where $z$ is the vertical displacement of the micro-beam, and the derivation operation is aimed at the variable $\tau$. Equation parameters $m$, $b$, $k_1$ and $k_3$ are effective lumped mass, damping coefficient, linear mechanical stiffness and cubic mechanical stiffness of the system, respectively. $G(\tau)$ is a random process can be specified into power law noise in the following discussion. The external driving force $F$ can be given by

$$F(z,\tau) = \frac{1}{2}\frac{C_0}{(d-z)^2}(V_b + V_{ac}\sin\Omega\tau)^2 - \frac{1}{2}\frac{C_0}{(d+z)^2}V_b^2, \tag{2}$$

where $C_0$ denotes the capacitance of the parallel-plate actuator with $z=0$, and $d$ means the initial gap width as well as $V_b$ the bias-voltage. $V_{ac}$ and $\Omega$ are the amplitude and frequency of the actuated alternating current voltage, respectively.

For the sake of convenience, it's sensible to introduce non-dimensional variables as follows:

$$t = \omega_0\tau, \; \omega = \frac{\Omega}{\omega_0}, \; x = \frac{z}{d}, \; \mu = \frac{b}{m\omega_0}, \; \delta = \frac{k_1}{m\omega_0^2}, \; \beta = \frac{k_3 d^2}{m\omega_0^2}, \; \gamma = \frac{C_0 V_b^2}{2m\omega_0^2 d^3}, \; A = 2\gamma\frac{V_{ac}}{V_b},$$

where $\omega_0 = \sqrt{\frac{k_1}{m}}$ can be seen as the natural frequency and the assumption that the amplitude of actuated AC voltage is much smaller than the bias voltage must be made when the MEMS resonator possesses high quality [16]. So the dimensionless equation of motion is acquired

$$\ddot{x} + \mu\dot{x} + \delta x + \beta x^3 = \gamma\left(\frac{1}{(1-x)^2} - \frac{1}{(1+x)^2}\right) + \frac{A}{(1-x)^2}\sin\omega t + g(t), \tag{3}$$

where the derivation operation is aimed at the variable $t$. Let $y = \dot{x}$, rewriting kinematic equation in the following form:

$$\begin{cases} \dot{x} = y \\ \dot{y} = -\delta x - \beta x^3 - \mu y + \gamma\left(\frac{1}{(1-x)^2} - \frac{1}{(1+x)^2}\right) + \frac{A}{(1-x)^2}\sin\omega t \end{cases}. \tag{4}$$

## 2.2 power law noise

The power law noise, also known as $1/f^\alpha$ noise, is defined according to the function form of its power spectral density $S(f)$. The measured spectral attenuation is approximatively $f^{-\alpha}$, where $f$ denotes the frequency and parameter $\alpha \in [0,2]$ is usually referred to as frequency exponent or attenuation index [23]. When the value of $\alpha$ is not demarcated strictly, the noise can be treated as $1/f$ noise in a broad sense, and the generalized noise can be pursued the trace in the drive system under nonequilibrium conditions. Within the scope of low frequency, the amplitude of power law noise is larger than other common noises, thus it receives much attention from investigators [24].

Power law noise has been observed in a myriad of scientific disciplines. White noise with $\alpha = 0$ is the most common, indeed, practically ubiquitous type of noise in the field of life. For example, when thermal noise generated by random motion of the electrons pass through the electrical resistance, it may cause the "snow" on radar and television screen if there is no input [25]. As for $1/f$ noise, it gains the beauteous appellation pink noise for that brown noise is also referred to as red noise and pink is neutralized between white and red [23]. Facts have demonstrated that $1/f$ fluctuations coincide with the variation rules of biological signals such as $\alpha$ brain waves when people in silent mode and heart rate cycle. In natural phenomena such as intensity changes in wind, light, average temperature on the surface of the earth and cosmic radiation, researchers have also observed $1/f$ fluctuation. But in the field of electronic science, $1/f$ fluctuations are seen as flick noise more often. The exploration and research of flick noise have a very long history. Johnson [26] successfully measured white noise spectrum and unexpected "flick noise" in the low frequency range. Through further experiments, Johnson discovered the power spectral density of current noise was inversely proportional to the frequency $f$, so he described it as $1/f$ noise and put forward the term "$1/f$ noise" for the first time. In addition to what mentioned above, Clarke and Voss found that both voice and music have $1/f$ spectra [27], and they even came up with an algorithm to compose a list of "fractal" music [28]. Gutenberg and Richter [29], the mainstays of the physics of earthquakes, they presented a law stating that $N(M)$, the number of earthquakes with magnitude large than $M$, is proportional to $10^{-bM}$, namely, $\log_{10} N(M) = A - bM$, where $b$ is the slope near to one. Musha [25], a physicist at the Tokyo Institute of Technology, discovered that it would show $1/f$ fluctuation when traffic flow past a certain site on a highway. Other $1/f^\alpha$ noise signals with $\alpha \neq 0$ or $1$ have been observed in statistical analyses in astronomy [30], electronic devices [31], DNA sequences [32], brain signals, biomedicine [33,34]，psychological mental states, and natural images. For example, Careri and Consolini [33] have observed the fluctuation of the electrical dipole moment of lysozome, an important enzyme whose spectral density is $1/f^\alpha$ spectrum, where $\alpha \approx 1.5$. In 1986, EXOSAT [30] monitored the X-ray lightcurve of the Seyfert Galaxy NGC 5506, and the slope of corresponding spectrum shows $1.8 \pm 0.3$. When $\alpha = 2$, the noise is usually referred to a series of colored noises with very strong correlation, in which case people will feel tedious and plodding. More detailed statements and introduction can be found in related articles.

Considering the actual environment, it's of great significance to provide pure power law signals for comparative analysis. In general, $1/f$ low frequency noise can be detected by instruments directly, but the noise obtained with this method may be contaminated by other types of noise. However, the sequences acquired by computer numerical simulation can evade this shortcoming. There are two main theories for generating $1/f$ low frequency noise, fractional Brownian movement mode [35] and Lorentz spectrum method [36]. In this study we choose the algorithm introduced by Miro Stoyanov et al. [23], and the core of the algorithm is the Z and fast Fourier transform with computational complexity of $O(n \log n)$. Compared with the two classical algorithms, the preponderance of the new digital model lies in allowing very precise and efficient computer generation of power law noise for any frequency exponent [37]. Figure 1 shows numerical simulation of power law noise series with various $\alpha$ and corresponding power spectrum densities.

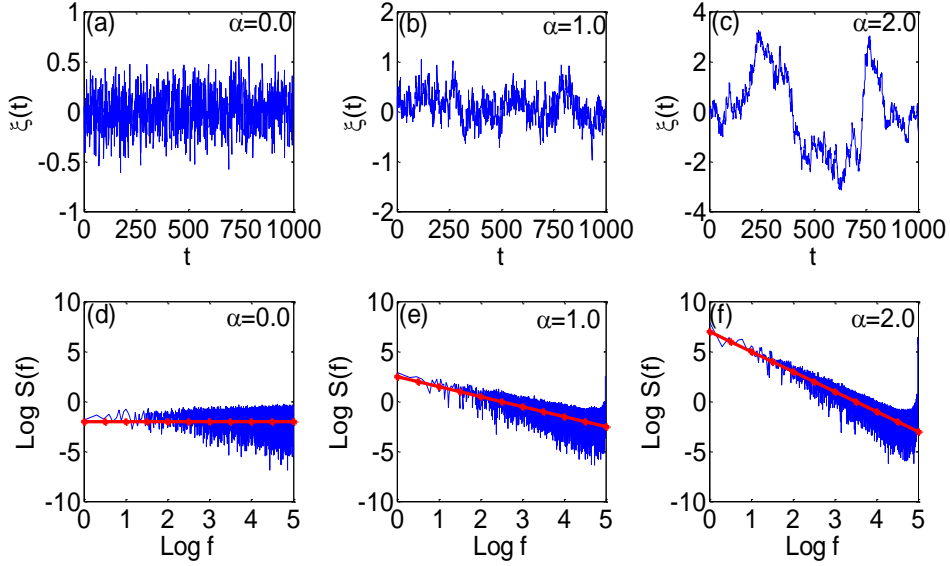

**Fig.1.** Simulation of discrete power law noise sequences with different $\alpha$ and corresponding power spectrum densities. Top row: noise sequence diagrams (a), (b), (c). Bottom row: Power spectrum density diagrams, (d), (e), (f). With the increment of $\alpha$, the simulation of noise sequence becomes smoother, indicating that the correlation of random process is getting stronger. For power spectrum density diagrams, the plot of log ($S(f)$) versus log (frequency) will be linear with slope, where the slope of straight line is frequency exponent.

## 3. Chaos prediction in MEMS resonators under power law noise

On the basis of the stochastic Melnikov method, the intensity threshold of power law noise will be deduced analytically in this section.

When taking the effect of noise intensity on the nonlinear dynamic system into account, the governing equation of motion can be given as

$$\begin{cases} \dot{x} = y \\ \dot{y} = -\delta x - \beta x^3 + \gamma\left(\dfrac{1}{(1-x)^2} - \dfrac{1}{(1+x)^2}\right) + \varepsilon\left(-\bar{\mu} y + \dfrac{\bar{A}}{(1-x)^2}\sin \omega t + \sigma \xi(t)\right) \end{cases}, \quad (5)$$

where $\bar{\mu} = \mu/\varepsilon$, $\bar{A} = A/\varepsilon$ ($\mu = O(\varepsilon)$, $A = O(\varepsilon)$), $\varepsilon$ is a small disturbance parameter, $\xi(t)$ represents power law noise and $\sigma$ is the noise intensity.

When $\varepsilon = 0$, Eq.(5) can be taken as an unperturbed system with the form:

$$\begin{cases} \dot{x} = y \\ \dot{y} = -\delta x - \beta x^3 + \gamma\left(\dfrac{1}{(1-x)^2} - \dfrac{1}{(1+x)^2}\right) \end{cases}. \quad (6)$$

System (6) is a Hamilton system and the corresponding Hamilton function is

$$H(x, y) = \frac{1}{2} y^2 + \frac{1}{2}\delta x^2 + \frac{1}{4}\beta x^4 - \gamma\left(\frac{1}{1-x} + \frac{1}{1+x}\right) + 2\gamma. \quad (7)$$

The form of the potential function is as follows when the potential energy is defined as zero at the point $x = 0$:

$$V(x) = \frac{1}{2}\delta x^2 + \frac{1}{4}\beta x^4 - \gamma\left(\frac{1}{1-x} + \frac{1}{1+x}\right) + 2\gamma. \quad (8)$$

It can be learned that, however, there is no accurate closed solution for the homoclinic trajectory of Eq.(5). The Melnikov method requires that the saddle points are connected by homoclinic orbit and the orbit can be expressed analytically. In order to utilize the method, it is optional to approximate it by expanding the unperturbed term $\left(\dfrac{1}{(1-x)^2} - \dfrac{1}{(1+x)^2}\right)$ to $\left(4x + 8x^3 + O(x^5)\right)$ in the form of Taylor series [16], and the homoclinic orbits will turn to be the typical frictionless Duffing equation. Based on the same thought, the term $\left(\dfrac{1}{(1-x)} + \dfrac{1}{(1+x)}\right)$ in the right side of Eq.(8) can be replaced by $\left(2 + 2x^2 + 2x^4 + O(x^5)\right)$. When put the replacement operation into effect, Eq.(8) will be rewritten as

$$V(x) = -\frac{1}{2}\kappa x^2 + \frac{1}{4}\lambda x^4, \quad (9)$$

where $\kappa = 4\gamma - \delta$ and $\lambda = \beta - 8\gamma$.

According to the signs of the parameters $\kappa$ and $\lambda$, the system can be classified into four categories. We only pay attention to the case where $\kappa > 0$ and $\lambda > 0$ $(\delta/4 < \gamma < \beta/8)$. Under these circumstances, one can get that the model is bistable, and at the same time, there are three equilibrium points, one unstable saddle point $P_0(0,0)$ and two stable center points $P^\pm(\pm x_p, 0)$.

When $\delta/4 < \gamma < \beta/8$, the homoclinic trajectory of Eq.(5) can be approximated by

$$\left(x_0(t), y_0(t)\right) = \left(\pm x_p \operatorname{sech}\left(\sqrt{\kappa} t\right), \mp \sqrt{\kappa} x_p \operatorname{sech}\left(\sqrt{\kappa} t\right) \tanh\left(\sqrt{\kappa} t\right)\right), \tag{10}$$

where $x_p = \Lambda\sqrt{\kappa}$, the positive root of Eq.(9), and $\Lambda = \sqrt{2/(\beta - 8\gamma)}$. Figure 2 shows the approximated homoclinic trajectory obtained with Eq.(10). Compared with the analytical solution of the model (5), the imitative effect is favorable, so the replacement operation is advisable.

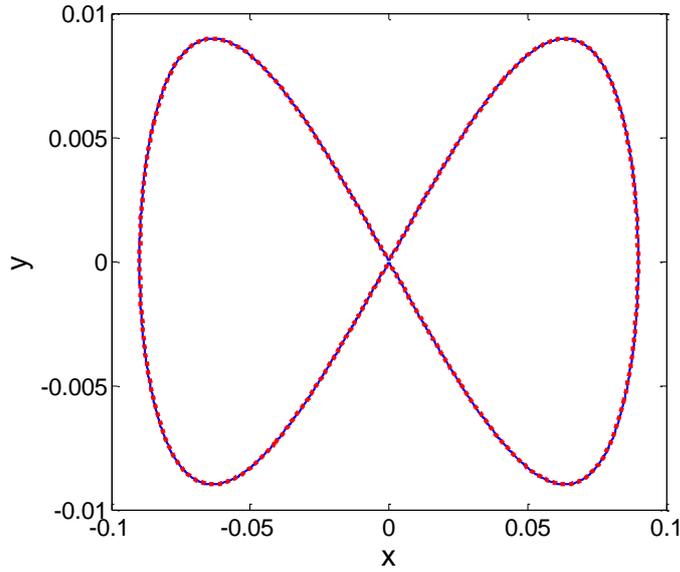

**Fig.2.** Comparison of numerical simulation (dashed line) and analytical result (solid line) of homoclinic trajectory in phase space. $\alpha = 1.0, \beta = 12.0$ and $\gamma = 0.26$.

By putting the method and formula presented by Wiggins [38] into use, the random Melnikov process can be obtained as

$$M(t_0) = \int_{-\infty}^{+\infty} -\bar{\mu} y_0^2(t) dt + \int_{-\infty}^{+\infty} \frac{\bar{A}}{(1-x_0(t))^2} \sin\omega(t+t_0) y_0(t) dt + \int_{-\infty}^{+\infty} \sigma y_0(t) \xi(t+t_0) dt \quad (11)$$
$$= I_d + I_p + I_s$$

where

$$I_d + I_p = \int_{-\infty}^{+\infty} -\bar{\mu} y_0^2(t) dt + \int_{-\infty}^{+\infty} \frac{\bar{A}}{(1-x_0(t))^2} \sin\omega(t+t_0) y_0(t) dt$$
$$= -\frac{2}{3} \bar{\mu} \sqrt{\kappa} x_p^2 + \bar{A} x_p \frac{2\pi\bar{\omega}}{\sqrt{1-x_p^2}} \cos(\bar{\omega} t_0) \frac{\sinh(-\bar{\omega}(\arccos(x_p)))}{\sinh(\bar{\omega}\pi)},$$

represents the deterministic part of the random Melnikov process due to damping and restoring force and the periodic excitation, and $I_s = \int_{-\infty}^{+\infty} \sigma y_0(t) \xi(t+t_0) dt$ denotes the random part with power law noise. The random Melnikov integral is an efficient measure of the random distance between stable and unstable manifolds. Once the distance is equal to zero, the transversal intersection will occur numerous times, which indicates the appearance of chaos [38].

Under the sense of mean value, one can get

$$E[(M(t_0))] = -\frac{2}{3} \bar{\mu} \sqrt{\kappa} x_p^2 + \bar{A} x_p \frac{2\pi\bar{\omega}}{\sqrt{1-x_p^2}} \cos(\bar{\omega} t_0) \frac{\sinh(-\bar{\omega}(\arccos(x_p)))}{\sinh(\bar{\omega}\pi)}, \quad (12)$$

where $E[\bullet]$ denotes ordinary expectation operator, and in this sense, chaos will never appear. So we consider the random Melnikov process (11) has simple zero points in the mean-square sense. When $\xi(t)$ is seen as an input of system (5), the impulse response function can be given as $h(t) = y_0(t)$ and the corresponding frequency response function is

$$H(\omega) = \int_{-\infty}^{+\infty} h(t)e^{-i\omega t} dt$$

$$= \pm \int_{-\infty}^{+\infty} i\sqrt{\kappa} x_p \operatorname{sech}(\sqrt{\kappa}t) \tanh(\sqrt{\kappa}t) \sin(\omega t) dt$$

$$= \pm i x_p \frac{\omega}{\sqrt{\kappa}} \pi \operatorname{sech}\left(\frac{\pi\omega}{2\sqrt{\kappa}}\right) \tag{13}$$

$$= \pm i \Lambda \pi \omega \operatorname{sech}\left(\frac{\pi\omega}{2\sqrt{\kappa}}\right)$$

where $i = \sqrt{-1}$. Thus, the variance of $I_s$, as the output of the system, can be obtained in frequency domain as

$$I_s^2 = E\left[\left(\sigma \int_{-\infty}^{+\infty} y_0(\tau)\xi(\tau+\tau_0)\right)^2\right]$$

$$= \sigma^2 \int_{-\infty}^{+\infty} |H(\omega)|^2 S(\omega) d\omega$$

$$= \sigma^2 \int_{-\infty}^{+\infty} \Lambda^2 \omega^2 \pi^2 \operatorname{sech}^2\left(\frac{\pi\omega}{2\sqrt{\kappa}}\right) \frac{1}{\omega^\alpha} d\omega \tag{14}$$

$$= \sigma^2 \Lambda^2 \pi^2 \int_{-\infty}^{+\infty} \operatorname{sech}^2\left(\frac{\pi\omega}{2\sqrt{\kappa}}\right) \omega^{2-\alpha} d\omega$$

From the standpoint of energy, the criterion for chaos is that the random Melnikov process has simple zero points in the mean-square sense. In other words, the mean-square formula is

$$I_d^2 + I_p^2 = I_s^2. \tag{15}$$

Then the possible condition for the occur of the chaos in system (5) can be expressed as

$$\left(-\frac{2}{3}\bar{\mu}\sqrt{\kappa}x_p^2\right)^2 + \left(\bar{A}x_p \frac{2\pi\bar{\omega}}{\sqrt{1-x_p^2}} \cos(\bar{\omega}t_0) \frac{\sinh\left(-\bar{\omega}(\arccos(x_p))\right)}{\sinh(\bar{\omega}\pi)}\right)^2$$

$$= \sigma^2 \Lambda^2 \pi^2 \int_{-\infty}^{+\infty} \operatorname{sech}^2\left(\pi\omega/2\sqrt{\kappa}\right) \omega^{2-\alpha} d\omega \tag{16}$$

The integral in the right side of Eq.(16) can be obtained with numerical methods. Eq.(16) predicts that power law noise can evoke chaotic response in dynamic systems and the chaotic motion is under the sense of Smale horseshoes mapping. In the following parts, we choose the fixed parameters as $\delta=1.0, \beta=12.0, \gamma=0.26, \varepsilon=0.1, \mu=0.01, V_b=3.8V, V_{ac}=0.13V, \omega=0.5$ unless otherwise indicated. The threshold of the noise intensity for the onset of chaos in system (5) is shown in Fig.3. The numerical results are indicated by discrete star points. Compared with the analytical results, sharp computational cost is primary weakness of the numerical method. As it is shown, the threshold of noise intensity decreases with the increasing of the frequency exponent $\alpha$. Furthermore, both the theoretical and numerical simulations show the same trend. According to the previous discussion, as the frequency exponent $\alpha$ increases, the correlation of noise sequence will enhance. This manifests that one can restrain the generation of chaos by strengthening the randomness of noise. From Fig.3, one can also acquire that when the frequency exponent is fixed, the system changes from order to chaos with the increasing of noise intensity, which can be easily attributed to the effect of noise. Meanwhile when fixing the noise intensity, the system also changes from order to chaos with the increasing of frequency exponent $\alpha$. Under the sense of Smale horseshoes mapping, chaos will only occur when noise intensity exceeds the threshold value. So, the injection of the noise can lead to chaotic behaviors, and a larger value of $\alpha$ or $\sigma$ will make this easier.

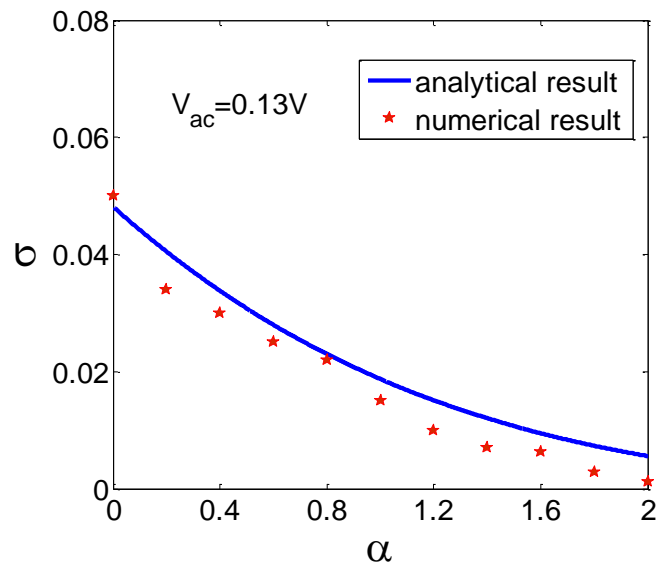

**Fig.3.** The threshold curve of the noise intensity $\sigma$ versus the frequency exponent $\alpha$. The numerical results are indicated by discrete star points, which fall region near the analytical curve.

## 4. Numerical simulations

Now the threshold of onset of chaos obtained with the stochastic Melnikov process in the sense of mean-square is only a possible criterion and the origin of chaotic dynamics cannot be guaranteed, while the horseshoe chaos is a transient phenomenon. So it is indispensable for us to use other numerical simulation methods to illustrate and certificate the analytical results.

### 4.1 Mean largest Lyapunov exponent

It is well known the fundamental feature of chaotic system is the hypersensitivity to initial conditions, namely the butterfly effect. When the initial values are given a negligible disturbance then the perturbed orbit will separate exponentially from the original one as time goes on. In theory, the long term behavior of dynamic system can't be controlled or predicted. Lyapunov exponent is a resultful measure describing the mean exponential rate of convergence or divergence between two adjacent orbits of phase space qualitatively and quantitatively. It is deemed as an indicator to determine whether chaotic behaviors exist in dynamic systems. If the largest Lyapunov exponent is less than zero, the system can be taken as quasi-periodic or periodic, i.e., steady state. On the contrary, chaotic systems with positive Lyapunov exponents may have strange attractors. For an n-dimensional dynamic system, one can get $n$ Lyapunov exponents by numerical calculation. The sum of overall Lyapunov exponents is equal to zero for the Hamilton system, but for the dissipative system, the sum is negative. No matter whether the system is dissipative or not, as long as the largest one is positive, it must exhibit chaotic behaviors.

With respect to the calculation of Lyapunov exponents, several methods are frequently selected. For example, the definition method, the orthogonal method, the Wolf's method and the small data set, etc. In this work, we utilize the Wolf's method [39] to calculate the Lyapunov exponents. The crux of the algorithm is to supervise the long-term development of an infinitesimal n-sphere of initial conditions. The largest Lyapunov exponent is defined as

$$L = \lim_{t \to \infty} \frac{1}{t} \ln \frac{p(t)}{p(0)}, \quad (17)$$

where $p(t)$ is the length of the principle axis. However, for stochastic systems, we employ the mean largest Lyapunov exponent for the sake of balancing the effect of random effects. Through calculating the largest Lyapunov exponents of all sample orbits, the mean largest Lyapunov exponent is

$$L = \frac{1}{N} \sum_{i=1}^{N} L_i, \quad (18)$$

where the $N$ denotes the number of orbits. In this paper, we choose $N = 100$ for research. In order to compare with the results of the stochastic Melnikov method, the corresponding numerical results of the chaotic thresholds are shown in Fig.3 by vanishing largest Lyapunov exponent. The two methods agree well with each other.

### 4.2 Phase diagrams and time histories

Further, we utilize two frequently used numerical methods, phase diagrams and time histories to underpin the conclusions obtained above. Phase diagrams contain all possible orbits in phase state, and one can determine whether there is chaotic motion in a specified system through the shape of phase trajectory. The time histories are utilized to verify the sensitive dependence on the initial conditions. The consequences are presented in Fig.4. As it is shown, when $\alpha = 0.0$, the shape of phase orbit changes from stable limit cycle to tanglesome line with the increasing of the noise intensity $\sigma$ from 0.03 to 0.1. In the meantime, the corresponding time history changes from organized periodic state to chaotic state. It is necessary to point out that transient time $t = 2400T$ is get rid of in numerical simulation, where T is the period of periodic excitation, so we can deem that Fig.4 represents a stable state. This is consistent with the result of the mean largest Lyapunov exponent method (the corresponding values $L \approx -0.0062, 0.0126$). Now, we take another circumstance into account. When $\sigma = 0.03$, the shape of phase trajectory and time history will go through the same transformation with adding the frequency exponent $\alpha$ from 0.0 to 1.2 (the corresponding mean largest Lyapunov exponent $L \approx -0.0062, 0.0258$).

As mentioned above, one can bring in a verdict whether there is chaotic motion in a given system by contrasting the tracks of time history under two closed initial conditions. As it is shown in Fig.4(b) that when the initial value is given a small scale disturbance, the time histories will have little difference and the system is still periodic. However, remarkable differences will be seen in Fig.4(d) and Fig.4(f) if we set the same perturbed value, which manifests a strong sensitive dependence on initial conditions. It is a significant fundamental feature of chaos.

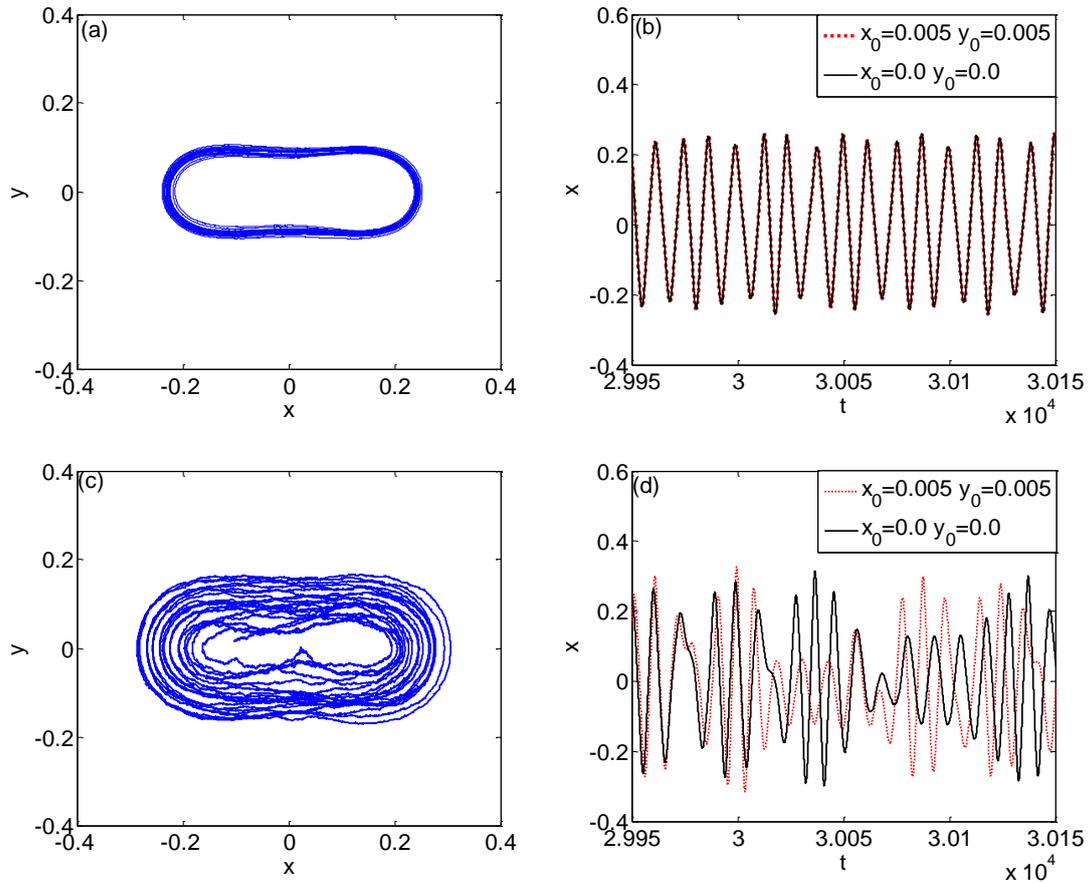

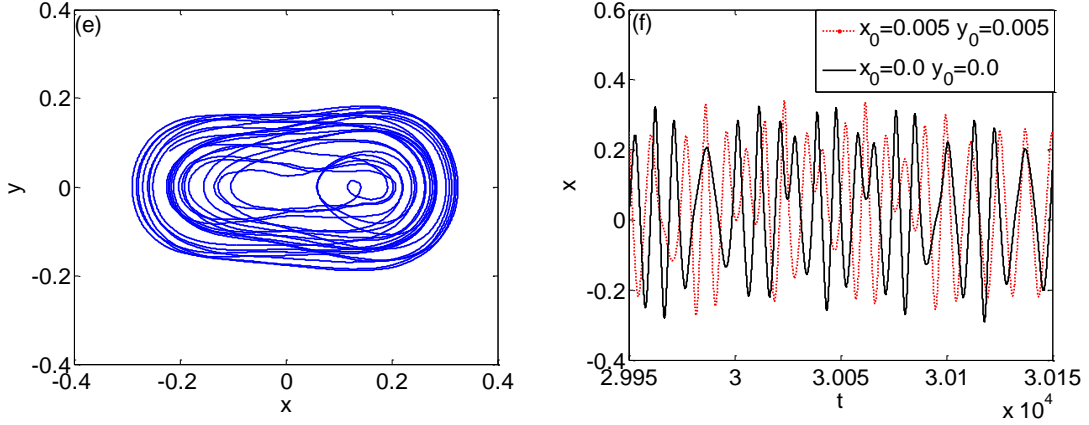

**Fig.4.** Phase diagrams and time histories for different initial conditions with $V_{ac}=0.13V$ in the system (6) respectively: ((a)(b)): $\sigma=0.03, \alpha=0.0$, ((c)(d)) $\sigma=0.1, \alpha=0.0$, ((e)(f)) $\sigma=0.03, \alpha=1.2$.

# 5. Chaos control of MEMS resonators with time-delay feedback method

In this section, a time-delay feedback control policy is applied to control the chaotic motion of MEMS resonator under power law noise.

As presented in the introduction, the existence of noise will have a significant effect on the properties and performance of MEMS resonators. So as to guarantee the normal operation of the system, it's necessary to exert chaos control strategies to lower the damage caused by noise. There are two main ideas of chaos control: one is that there is no specific control object, but the chaos is eliminated by lower the Lyapunov exponent without considering the final state of the system; the other one is to transform the chaotic motion of the controlled dynamic system into prespecified periodic dynamic behavior by exerting distinct exterior controls. As for chaos control methods, the classical ones include OGY method [40], OPF method [41], adaptive control method [42] and time-delay continuous feedback control [43]. Based on the OGY control strategy, a delayed feedback control method is proposed by Pyragas [44]. The core of the method is to directly extract the output signal of the system, and then feed back to the chaotic system as the control signal after certain time delay, without knowing the predetermined orbit in advance.

Let us add a time-delay feedback controller to the MEMS resonator model, and the motion equation can be expressed as

$$\begin{cases} \dot{x} = y \\ \dot{y} = -\delta x - \beta x^3 - \mu y + \gamma \left( \dfrac{1}{(1-x)^2} - \dfrac{1}{(1+x)^2} \right) + \dfrac{A}{(1-x)^2} \sin \omega t + g(t) + f(t) \end{cases}, \quad (19)$$

where $f(t)$ is the time-delay feedback controller and $f(t) = k(y(t - t_d) - y(t))$, $k = \varepsilon \bar{k}$ is the feedback gain factor and $t_d = 2\pi / \omega$ means the corresponding delay time. Similar to the previous mode, Eq.(19) can be showed in the form

$$\begin{cases} \dot{x} = y \\ \dot{y} = -\delta x - \beta x^3 + \gamma \left( \dfrac{1}{(1-x)^2} - \dfrac{1}{(1+x)^2} \right) + \varepsilon \left( -\bar{\mu} y + \dfrac{\bar{A}}{(1-x)^2} \sin \omega t + \sigma \xi(t) + \bar{f}(t) \right) \end{cases}, \quad (20)$$

where $f(t) = \varepsilon \bar{f}(t)$.

The stochastic Melnikov process from Eq.(20) can be expressed as

$$M(t_0) = \int_{-\infty}^{+\infty} y_0 \left[ -\bar{\mu} y_0 + \dfrac{\bar{A}}{(1-x_0)^2} \sin \omega(t + t_0) + \sigma \xi(t + t_0) + \bar{k}(y_{td} - y_0) \right] dt \\ = I_d + I_p + I_s + I_t \quad (21)$$

where

$$I_t = \int_{-\infty}^{+\infty} \bar{k} y_0 (y_{td} - y_0) dt = \bar{k} \kappa x_p^2 I_t - \dfrac{2\bar{k}}{3} \sqrt{\kappa} x_p^2. \quad (22)$$

From Eq.(21), the probable criterion for suppressing chaos in the mean-square sense is

$$I_d^2 + I_t^2 = I_p^2 + I_s^2. \quad (23)$$

Now it is essential to verify the effectiveness of the time-delay feedback controller. Here we set $\sigma = 0.1$, and the system is chaotic in the entire parameter domain of $\alpha$. Corresponding to the case indicated in Fig.4(c)(d), where the chaos control algorithm is not imposed, namely $k = 0$, the Lyapunov exponent image is showed in Fig.5(a). We can discern that the system is in a chaotic state when $\sigma = 0.1$. While the chaos controller exerted, the relationship between Lyapunov exponent and feedback gain factor $k$ is provided in Fig.5(b). With the increased intensity of feedback control, the Lyapunov exponent decreases monotonically. When the control element is greater than 0.06, the dynamic system will change from chaotic state to ordered state, which demonstrates the effectiveness of the chaos control algorithm. Echoing with the Fig.5(b), the phase diagrams and time histories are shown in Fig.6.

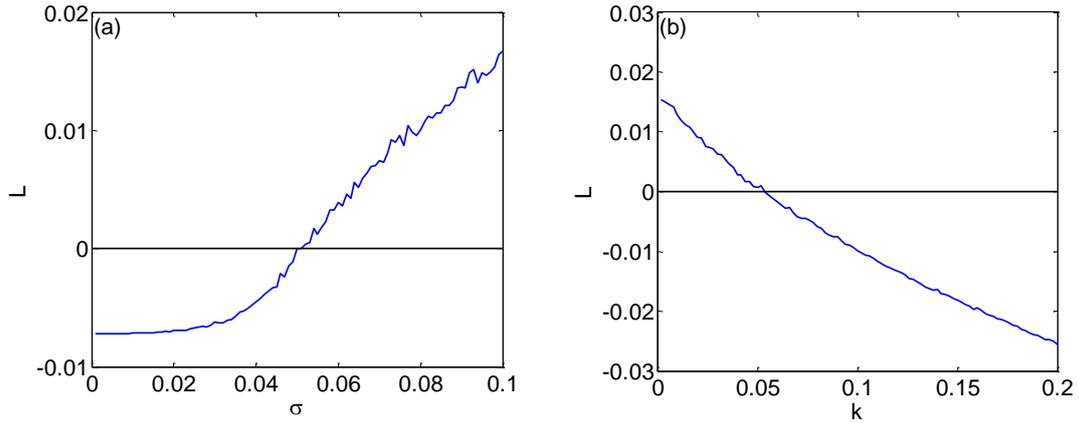

**Fig.5.** The Lyapunov exponent curve for different states of the system. (a): no chaos control item added. (b): time-delay feedback controller imposed.

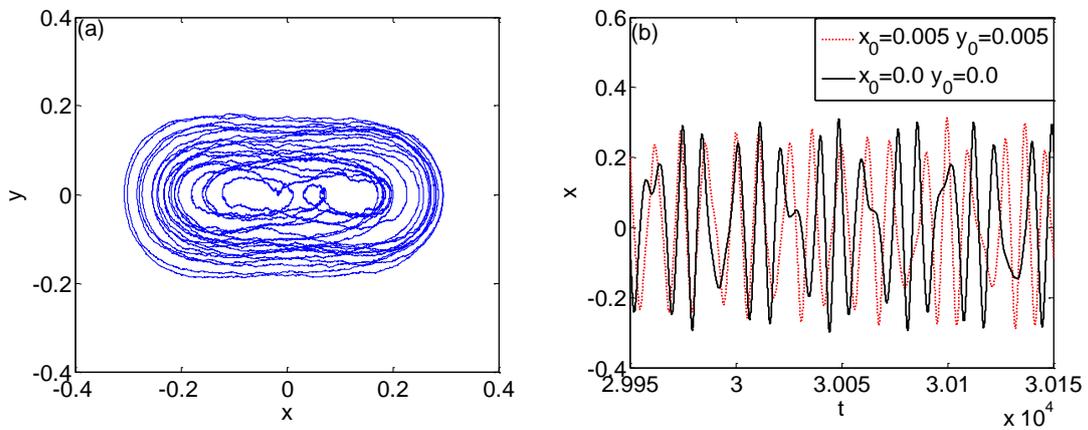

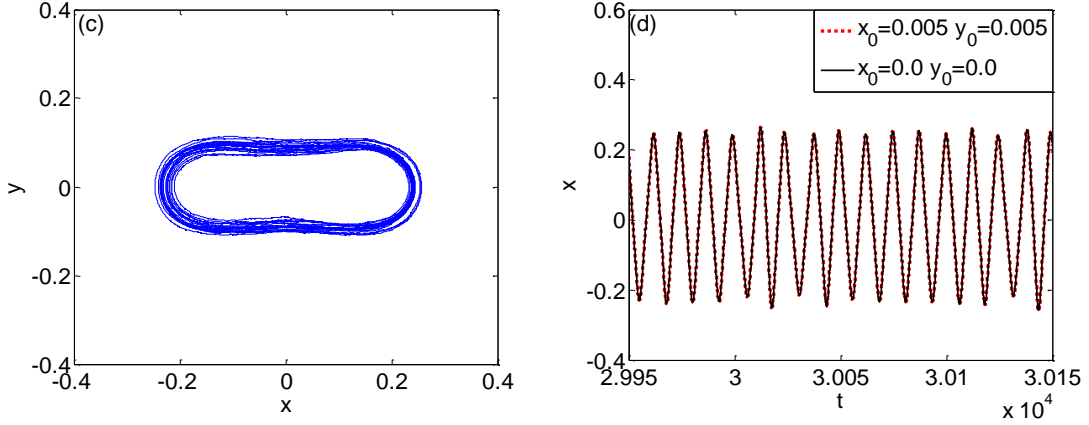

**Fig.6.** Phase diagrams and time histories for different initial conditions with $V_{ac} = 0.13V$ in the system (21) respectively: ((a)(b)): $\sigma = 0.1, \alpha = 0.0, k = 0.0$, ((c)(d)): $\sigma = 0.1, \alpha = 0.0, k = 0.2$.

# 6. Conclusion

In this work, chaos and chaos control in the electrostatically actuated MEMS resonators are reported analytically and numerically. Other than electrostatic force, the random disturbance, namely power law noise, is also exerted on the clamped-clamped microbeam. By simple Taylor expansion, we acquire the approximated homoclinic orbit of the dynamic system in the form of frictionless duffing equation. Compared with the analytical solution of the model, the imitative effect is favorable, so the replacement operation is recommended. So as to gain a more precise solution, high order terms of approximation must be considered when solving the homoclinic trajectory. On the basis of approximation concept and the mathematical model, the kinetic equation is applicable for the universal stochastic Melnikov method. Through rigorous derivation and analysis, we realize that the power law noise can really induce stochastic chaotic motion and the threshold of noise intensity decreases with the increasing of the frequency exponent of power law noise in entire parameter field. We demonstrate that when the noise intensity exceeds a threshold value, stochastic chaos will occur in the system. The analytical results are certified by numerical simulations in light of the mean largest Lyapunov exponent, phase diagram and time history. Based on the previous experimental results, a time-delay feedback control algorithm is proposed for controlling chaotic motion. Similarly, the criterion of chaos control is derived by the stochastic Melnikov method. The effectiveness of this controller is also certificated by the above numerical method. In practice, one can induce or restrain the occurrence of chaos in the MEMS resonators by altering noise environment and exerting the time-delay feedback control.

## Acknowledgements

The work is supported by the National Natural Science Foundation of China (Grant No. 11672231) and the NSF of Shaanxi Province (Grant No. 2016JM1010).